\input amstex
\documentstyle{amsppt}
\magnification=\magstep1
 \hsize 13cm \vsize 18.35cm \pageno=1
\loadbold \loadmsam
    \loadmsbm
    \UseAMSsymbols
\topmatter
\NoRunningHeads
\title
Euler numbers and polynomials associated with zeta functions
\endtitle
\author
  Taekyun Kim
\endauthor
 \keywords Bernoulli numbers,twisted Bernoulli numbers, zeta function,
Euler numbers, twisted Euler numbers, $q$-Volkenborn integral,
Dedekind sum, $p$-adic invariant integral
\endkeywords

\abstract For $s\in\Bbb C$, Euler zeta function and Hurwitz's type
Euler zeta function are defined by
$$\zeta_{E}(s)=2\sum_{n=1}^{\infty}\frac{(-1)^{n}}{n^s}, \text{ and } \zeta_{E}(s,x)=2\sum_{n=0}^{\infty}
\frac{(-1)^n}{(n+x)^s}, \text{ see [ 20, 21, 22, 31, 15 ]}.$$ Thus,
we note that Euler zeta functions are entire functions in whole
complex $s$-plane and these zeta functions have the values of Euler
numbers or Euler polynomials at negative integers. That is,
$$ \zeta_{E}(-k)=E_k^*, \text{ and  } \zeta_{E}(-k,x)=E_k^*(x), \text{ cf. [ 12, 20, 21, 22, 31] }. $$

In this paper, we give some interesting identities between Euler
numbers and zeta functions. Finally we will give the new values of
Euler zeta function at positive even integers.

\endabstract
\thanks  2000 AMS Subject Classification: 11B68, 11S80
\newline
\endthanks
\endtopmatter

\document

{\bf\centerline {\S 1. Introduction}}
 \vskip 20pt

Throughout  this paper $\Bbb Z$, $\Bbb Q$, $\Bbb C$,  $\Bbb Z_p$,
$\Bbb Q_p$ and $\Bbb C_p$ will respectively denote the ring of
rational integers,  the field of rational numbers, the filed of
complex numbers, the ring $p$-adic rational integers, the field of
$p$-adic rational numbers, and the completion of the algebraic
closure of $\Bbb Q_p$. Let $v_p$ be the normalized exponential
valuation of $\Bbb C_p$ such that $|p|_p=p^{-v_p(p)}=p^{-1}$. If
$q\in\Bbb C_p$, we normally assume $|q-1|_p<1$. We use the notation
$$[x]_q=\frac{1-q^x}{1-q}, \text{ and  }
[x]_{-q}=\frac{1-(-q)^x}{1+q}.$$ Hence, $\lim_{q\rightarrow 1}[x]_q
=1, $ for any $x$ with $|x|_p\leq 1$ in the present $p$-adic case.

Let $p$ be a fixed odd prime. For $d(= odd)$ a fixed positive
integer with $(p,d)=1$, let
$$\split
& X=X_d = \lim_{\overleftarrow{N} } \Bbb Z/ dp^N \Bbb Z , \ X_1 =
\Bbb Z_p , \cr  & X^\ast = \underset {{0<a<d p}\atop {(a,p)=1}}\to
{\cup} (a+ dp \Bbb Z_p ), \cr & a+d p^N \Bbb Z_p =\{ x\in X | x
\equiv a \pmod{dp^n}\},\endsplit$$ where $a\in \Bbb Z$ lies in
$0\leq a < d p^N$.

In this paper we prove  that $$\mu_{-q}(a+dp^N \Bbb
Z_p)=(1+q)\frac{(-1)^aq^a}{1+q^{dp^N}}
=\frac{(-q)^a}{[dp^N]_{-q}},$$ is distribution on $X$ for $q\in \Bbb
C_p$ with $|1-q|_p< 1. $ This distribution yields an integral as
follows:
$$I_{-q}(f)=\int_{\Bbb Z_p} f(x)d\mu_{-q}(x)=\lim_{N\rightarrow
\infty}\frac{1}{[p^N]_{-q}}\sum_{x=0}^{p^N-1}f(x)(-q)^x, \text{ for
$f\in UD(\Bbb Z_p)$ },
$$
which has a sense as we see readily that the limit is convergent
(see [14]). Let $q=1$. Then we have the fermionic $p$-adic integral
on $\Bbb Z_p$ as follows:
$$I_{-1}=\int_{\Bbb
Z_p}f(x)d\mu_{-1}(x)=\lim_{N\rightarrow
\infty}\sum_{x=0}^{p^N-1}f(x)(-1)^x, \text{ cf.[1, 5, 20, 21, 14].}
$$
 For any positive integer $N,$ we set
$$\mu_q (a+lp^N \Bbb Z_p)=\frac{q^a }{[lp^N ]_q}, \text{ cf. [1, 5, 6, 7, 10, 11-18, 27-31]},$$
and this can be extended to a distribution on $X$. This distribution
yields $p$-adic bosonic $q$-integral as follows (see [13, 31]):
$$ I_q (f)=\int_{\Bbb Z_p}  f(x) \,d\mu_q(x)=\int_X  f(x) \,d \mu_q (x) ,$$
where $f \in UD(\Bbb Z_p )=\text{ the space of uniformly
differentiable function  on } \Bbb Z_p $ with values in $\Bbb C_p$,
cf.[1, 13, 27, 28, 29, 30, 31]. In view of notation, $I_{-1}$ can be
written symbolically as $I_{-1}(f)=\lim_{q\rightarrow -1} I_{q}(f).$
If we take $f(x)=q^{-x}[x]_q^n$, then we can derive the
$q$-extension of Bernoulli numbers and polynomials from $p$-adic
$q$-integrals on $\Bbb Z_p$ as follows:
$$\beta_{n,q}=\int_{\Bbb Z_p}q^{-x}[x]_q^nd\mu_q(x), \text{ and } \beta_{n,q}(x)=\int_{\Bbb Z_p}q^{-y}[y+x]_qd\mu_{q}(y)
, \text{ cf.[ 13, 31]}. $$ Thus, we note that
$$\beta_{0,q}=\frac{q-1}{\log q}, \text{ and } \beta_{m,q}=\frac{1}{(q-1)^m}\sum_{i=0}^m\binom{m}{i}\frac{i}{[i]_q},
\text{ cf.[ 13, 17, 31]. }$$ In complex plane, the ordinary
Bernoulli numbers are a sequence of signed rational numbers that can
be defined by the identity

 $$\frac{t}{e^t-1}=\sum_{n=0}^{\infty}B_n \frac{t^n}{n!}, \text{
 $|t|<2\pi$ }, \text{ cf.[1-33]}. $$

These numbers arise in the series expansions of trigonometric
functions, and are extremely important in number theory and
analysis. From the generating function of Bernoulli numbers, we note
that $B_0=1, B_1=-\frac{1}{2}, B_2=\frac{1}{6}, B_4=-\frac{1}{30},
B_6=\frac{1}{42}, B_8=-\frac{1}{30}, B_{10}=\frac{5}{66},
B_{12}=-\frac{691}{2730}, B_{14}=\frac{7}{6},
B_{16}=-\frac{3617}{510}, B_{18}=\frac{43867}{798},
B_{20}=-\frac{174611}{330}, \cdots,$ and $B_{2k+1}=0$ for $k\in\Bbb
N$. It is well known that Riemann zeta function is defined by
$$ \zeta(s)=\sum_{n=1}^{\infty}\frac{1}{n^s}, \text{ for $s\in\Bbb
C$ }.$$ We also note that Riemann zeta function is closely related
to Bernoulli numbers at positive integer or negative integer in
complex plane. Riemann did develop the theory of analytic
continuation needed to rigorously define $\zeta(s)$ for all $s\in
\Bbb C-\{0\}$. From this zeta function, he derived the following
formula, cf.[1-33]:
$$\zeta(-n)=-\frac{B_n+1}{n+1}, \text{ $n \in \Bbb N= \{1,2,3,\cdots \}$.} $$
Thus, we note that $\zeta(-n)=0$ if $n$ is even integer and greater
than 0. These are called the trivial zeros of the zeta function.
 In 1859, starting with Euler's factorization of the zeta
function
$$\zeta(s)=\prod_{p: prime}\frac{1}{1-p^{-s}},$$
he derived an explicit formula for the prime numbers in terms of
zeros of the zeta function. He also posed the Riemann Hypothesis: if
$\zeta(z)=0,$ then either $z$ is a trivial zero or $z$ lies on the
critical line $Re(z)=\frac{1}{2},$ cf.[20-33]. It is well known that
$$\frac{\sin z}{z}
=\left(1-\frac{z^2}{\pi^2}\right)\left(1-\frac{z^2}{(2\pi)^2}\right)\left(1-\frac{z^3}{(3\pi)^2}\right)
\cdots.$$ Thus,  $$1-z\cot
z=2\sum_{m=1}^{\infty}\frac{\zeta(2m)}{\pi^{2m}}z^{2m}, \text{
cf.[12, 20-33]}.$$ From this, we can derive the following famous
formula:

\proclaim{Lemma} For $n\in\Bbb N$, we have
$$\zeta(2n)=\sum_{k=1}^{\infty}\frac{1}{k^{2n}}=\frac{(-1)^{n-1}(2\pi)^{2n}}{2(2n)!}B_{2n},
\text{ for $n\in\Bbb N $} .$$
\endproclaim
 Since
$$ z \cot z=\frac{2iz}{e^{2iz}-1}+iz=1+\sum_{k=1}^{\infty}\frac{(-1)^k2^{2k}B_{2k}}{(2k)!}z^{2k}.$$
However, it is not known the values of $\zeta(2k+1)$ for $ k\in \Bbb
N$. In the case of $k=1$, Apery proved that $\zeta(3)$ is irrational
number. The constants $E_k^*$ in the Taylor series expansion
$$\frac{2}{e^t+1}=\sum_{n=0}^{\infty}E_n^{*}\frac{t^n}{n!}, \text{
where $|t|<\pi,$ cf.[5, 12, 20, 21, 22],} \tag1$$ are known as the
first kind Euler numbers. From the generating function of the first
kind Euler numbers, we note that
$$ E_0^*=1, \text{ and } E_n^*=-\sum_{l=0}^n\binom{n}{l}E_l^*, \text{
for $ n\in \Bbb N$ } .$$
 The first few are $1, -\frac{1}{2}, 0, \frac{1}{4},
\cdots,$ and $E_{2k}^{*}=0$ for $k=1,2,\cdots.$ The Euler
polynomials are also defined by
$$\frac{2}{e^t+1}e^{xt}=\sum_{n=0}^{\infty}E_n^{*}(x)\frac{t^n}{n!}=\sum_{n=0}^{\infty}\left(\sum_{k=0}^n
\binom{n}{k}E_k^* x^{n-k}\right)\frac{t^n}{n!}.$$ For $s\in\Bbb C$,
Euler zeta function and Hurwitz's type Euler zeta function are
defined by
$$\zeta_{E}(s)=2\sum_{n=1}^{\infty}\frac{(-1)^{n}}{n^s}, \text{ and } \zeta_{E}(s,x)=2\sum_{n=0}^{\infty}
\frac{(-1)^n}{(n+x)^s}, \text{ cf.[1, 11, 12, 20, 21, 22]}.$$ Thus,
we note that Euler zeta functions are entire functions in whole
complex $s$-plane and these zeta functions have the values of Euler
numbers or Euler polynomials at negative integers. That is,
$$ \zeta_{E}(-k)=E_k^*, \text{ and  } \zeta_{E}(-k,x)=E_k^*(x), \text{ cf.[1, 11, 12, 20, 21, 22] }. $$

In this paper, we give some interesting identities between Euler
numbers and zeta functions. Finally we will give the values of
Euler zeta function at positive even integers. \vskip 40pt

{\bf\centerline {\S 2. Preliminaries/Euler numbers associated with
$p$-adic fermionic integrals }}\vskip 10pt

Let $f_1(x)$ be translation with $f_1(x)=f(x+1)$. Then we have
$$ I_{-1}(f_{1})=I_{-1}(f)+2f(0).$$
If we take $f(x)=e^{(x+y)t}$, then we can derive the first kind
Euler  polynomials from the integral equation of $I_{-1}(f)$ as
follows:
$$  \int_{\Bbb Z_p}e^{(x+y)t}d\mu_{-1}(y)=e^{xt}\frac{2}{e^{t}+1}=\sum_{n=0}^{\infty}\frac{E_n^*(x)t^n}{n!}.$$
That is,
$$\int_{\Bbb Z_p}y^nd\mu_{-1}(y)=E_n^*, \text{ and} \int_{\Bbb Z_p}(x+y)^nd\mu_{-1}(y)=E_{n}^{*}(x).$$
For $n\in\Bbb N$,  we have the following integral equation:
$$\int_{\Bbb
Z_p}f(x+n)d\mu_{-1}(x)=(-1)^{n}\int_{\Bbb
Z_p}f(x)d\mu_{-1}(x)+2\sum_{l=0}^{n-1}(-1)^lf(l).$$ From this we
note that
$$\aligned
&E_{k}^{*}(n)-E_k^{*}=2\sum_{l=0}^{n-1}(-1)^{l-1}l^k, \text{ if
$n\equiv 0 (\mod 2)$, }\\
&E_k^{*}(n)+E_k^*=2\sum_{l=0}^{n-1}(-1)^ll^k, \text{ if $n \equiv
1 (\mod 2)$.}
\endaligned$$

Let $f(x)=\sin ax$ (or $f(x)=\cos ax $). By using the fermionic
$p$-adic $q$-integral on $\Bbb Z_p$, we see that
$$\aligned
0&=\int_{\Bbb Z_p} \sin ax d\mu_{-1}(x)+\int_{\Bbb Z_p}\sin ax
d\mu_{-1}(x)\\
&=(\cos a +1)\int_{\Bbb Z_p}\sin ax d\mu_{-1}(x)+\sin a\int_{\Bbb
Z_p}\cos ax d\mu_{-1}(x), \text{ see [15]}.
\endaligned$$
and
$$ 2=(\cos a +1)\int_{\Bbb Z_p} \cos ax d\mu_{-1}(x)-\sin
a\int_{\Bbb Z_p} \sin ax d\mu_{-1}(x). $$ Thus, we obtain
$$\int_{\Bbb Z_p}\cos ax d\mu_{-1}(x)=1, \int_{\Bbb Z_p}\sin ax
d\mu_{-1}(x)=-\frac{ \sin a}{\cos a+1}, \text{ see [15]}. $$ From
this we note that
$$\tan \frac{a}{2}=-\int_{\Bbb Z_p}\sin ax
d\mu_{-1}(x)=\sum_{n=0}^{\infty}\frac{(-1)^{n+1}a^{2n+1}}{(2n+1)!}E_{2n+1}^*.
$$
By the same motivation, we can also observe that
$$\frac{a}{2}\cot \frac{a}{2}=\int_{\Bbb Z_p} \cos ax
d\mu_{1}(x)=\sum_{n=0}^{\infty}\frac{(-1)^nB_{2n}}{(2n)!}a^{2n},
\text{ see [15].}
$$
These formulae are also treated in the Section 3.

Let $f(x)=e^{t(2x+1)}$. Then we can derive the generating function
of the second kind Euler numbers from fermionic $p$-adic integral
equation as follows:
$$\int_{\Bbb
Z_p}e^{t(2x+1)}d\mu_{-1}(x)=\frac{2}{e^{t}+e^{-t}}=\frac{1}{\cosh
t}=\sum_{n=0}^{\infty}E_n \frac{t^n}{n!}.$$ Thus, we have
$$(E+1)^n+(E-1)^n=2\delta_{0,n}, $$
where we have used the symbolic notation $E_n$ for $E^n$. The
first few are $E_0=1, E_1=0, E_2=-1, E_3=0, E_4=5, \cdots$,
$E_{2k+1}=0$ for $k\in\Bbb N$. In particular,
$$E_{2n}=-\sum_{k=0}^{n-1}\binom{2n}{2k}E_{2k}.$$
In the recent, Simsek, Ozden, Cangul, Cenkci, Kurt, etc have studied
the various extensions of the first kind Euler numbers by using
fernionic $p$-adic invariant $q$-integral on $\Bbb Z_p$, see [ 1, 5,
20, 21, 22, 27, 31]. It seems to be also interesting to study  the
$q$-extensions of the second kind Euler numbers due to Simsek et al(
see [20, 21, 27]).

\vskip 20pt

{\bf\centerline {\S 3. Some relationships between Euler numbers and
zeta functions }}\vskip 10pt

In this section we also consider Bernoulli and the second Euler
numbers in complex plane.  The second kind Euler numbers $E_k$ are
defined by the following expansion:
$$  sech x=\frac{1}{\cosh x}=\frac{2e^x}{e^{2x}+1}=\sum_{k=0}^{\infty}E_k\frac{t^k}{k!},
\text{ for $|t|<\frac{\pi}{2}$, cf.[12]}.\tag2$$ From (1) and (2),
we can derive the following equation:
$$E_k=\sum_{l=0}^k\binom{k}{l}2^lE_l^*,
\text{ where $\binom{k}{l}$ is binomial coefficient}. \tag3$$ By (3)
and (1), we easily see that $E_0=1, E_1=0, E_2=-1, E_3=0, E_4=5,
E_6=61, \cdots, $ and $E_{2k+1}=0$ for $k=1,2,3, \cdots.$ As Euler
formula, it is well known that
$$e^{ix}=\cos x+ i\sin x , \text{ where $i=(-1)^{\frac{1}{2}}$}.\tag4$$
From (4), we note that $\cos x=\frac{e^{ix}+e^{-ix}}{2}. $ Thus, we
have
$$\aligned
\sec x&= \frac{2}{e^{ix}+e^{-ix}}=sech
(ix)=\sum_{n=0}^{\infty}\frac{i^nE_n}{n!}x^n\\
&=\sum_{n=0}^{\infty}\frac{(-1)^nE_{2n}}{(2n)!}x^{2n}
+i\sum_{n=0}^{\infty}\frac{(-1)^nE_{2n+1}}{(2n+1)!}x^{2n+1}
=\sum_{n=0}^{\infty}\frac{(-1)^nE_{2n}}{(2n)!}x^{2n}.
\endaligned\tag5$$
From (5), we derive
$$x\sec
x=\sum_{n=0}^{\infty}\frac{(-1)^nE_{2n}}{(2n)!}x^{2n+1}, \text{
 for $|x|<\frac{\pi}{2}$}.\tag6$$
The Fourier series of an odd function on the interval $(-p, p)$ is
the sine series:
$$f(x)=\sum_{n=1}^{\infty}b_n \sin (\frac{n\pi x}{p}), \tag7$$
where
$$b_n=\frac{2}{p}\int_{0}^p f(x)\sin (\frac{n\pi x}{p}) dx.\tag8$$
Let us consider $f(x)=\sin ax$ on $[-\pi, \pi ]$. From (7) and (8),
we note that
$$ \sin ax =\sum_{n=1}^{\infty} b_n \sin nx ,\tag9$$
where
$$\aligned
b_n&= \frac{2}{\pi}\int_0^{\pi}\sin ax \sin nx
dx=\frac{2}{\pi}\int_0^{\pi}\left[ \frac{\cos
(n-a)x-\cos(n+a)x}{2}\right]dx\\
&=\frac{1}{\pi}\left[\frac{\sin (n-a)x}{n-a}-\frac{\sin
(n+a)x}{n+a}\right]_{0}^{\pi}=(-1)^{n-1}\frac{2}{\pi}\sin a\pi
\left(\frac{n}{n^2-a^2}\right).
\endaligned\tag10$$
In (9), if we take $x=\frac{\pi}{2}$, then we have
$$\aligned
&\sin \frac{\pi a}{2}=\sum_{n=1}^{\infty}b_{2n-1}(-1)^{n-1}\\
&=\frac{2}{\pi}\sin
a\pi\sum_{n=1}^{\infty}(-1)^{n-1}\frac{2n-1}{(2n-1)^2-a^2}
=\frac{2}{\pi}\sin a\pi\sum_{n=1}^{\infty}\frac{(2n-1)(-1)^{n-1}}{(2n-1)^2\left(1-(\frac{a}{2n-1})^2 \right)}\\
&=\frac{2}{\pi} \sin
a\pi\sum_{n=1}^{\infty}\frac{(-1)^{n-1}}{2n-1}\sum_{k=0}^{\infty}\frac{a^{2k}}{(2n-1)^{2k}}=\frac{2}{\pi}\sin
a\pi
\sum_{k=0}^{\infty}\left(\sum_{n=1}^{\infty}\frac{(-1)^{n-1}}{(2n-1)^{2k+1}}\right)a^{2k}
\endaligned\tag11$$
From (11), we note that
$$\frac{\pi a}{2}\sec (\frac{\pi
a}{2})=\sum_{k=0}^{\infty}\left(2\sum_{n=1}^{\infty}\frac{(-1)^{n-1}}{(2n-1)^{2k+1}}\right)a^{2k+1}.
\tag12$$ In (6), it is easy to see that
$$\frac{\pi a}{2}\sec(\frac{\pi
a}{2})=\sum_{n=0}^{\infty}\frac{(-1)^nE_{2n}}{(2n)!}\left(\frac{\pi}{2}\right)^{2n+1}a^{2n+1}.
\tag13$$ By (12) and (13), we obtain the following:

\proclaim{ Theorem 1} For $n\in \Bbb N$, we have

$$\sum_{k=1}^{\infty} \frac{(-1)^{k-1}}{(2k-1)^{2n+1}}
=\sum_{k=0}^{\infty} \frac{(-1)^{k}}{(2k+1)^{2n+1}}
=(-1)^{n}\frac{E_{2n}}{2(2n)!}\left(\frac{\pi}{2}\right)^{2n+1}.
\tag14$$
\endproclaim
It is easy to see that
$$\aligned\sum_{n=1}^{\infty}\frac{(-1)^n}{(2n+1)^{2k+1}}
&=2\sum_{n=1}^{\infty}\frac{1}{(4n-3)^{2k+1}}+\sum_{n=1}^{\infty}\frac{1}{(2n)^{2k+1}}
-\sum_{n=1}^{\infty}\frac{1}{n^{2k+1}}-1\\
&=\frac{1}{2^{4k+1}}\zeta(2k+1,
\frac{1}{4})-\frac{2^{2k+1}-1}{2^{2k+1}}\zeta(2k+1)-1.
\endaligned\tag15$$
By (14) and (15), we obtain the following:

\proclaim{ Corollary 2} For $n\in\Bbb N$, we have
$$\zeta(2n+1, \frac{1}{4})+2^{2n}(1-2^{2n+1})\zeta(2n+1)=(-1)^n
\frac{E_{2n}}{2(2n)!}\pi^{2n+1}2^{2n}.$$
\endproclaim

By simple calculation, we easily see that
$$i \tan
x=\frac{e^{ix}-e^{-ix}}{e^{ix}+e^{-ix}}=1-\frac{2}{e^{2ix}-1}+\frac{4}{e^{4ix}-1}.$$
Thus, we have
$$\aligned
x\tan x &=-xi +\frac{2xi}{e^{2xi}-1}-\frac{4xi}{e^{4xi}-1}\\
&=\sum_{n=1}^{\infty}\frac{(-1)^nB_{2n}4^n (1-4^n)}{(2n)!}x^{2n}.
\endaligned\tag16$$
From (16), we can easily derive
$$\tan x
=\sum_{n=0}^{\infty}\frac{(-1)^{n+1}4^{n+1}(1-4^{n+1})B_{2n+2}}{(2n+2)!}x^{2n+1}.
\tag17$$ By (4), we also see that
$$i \tan x
=1-\frac{2}{e^{2ix}+1}=i\sum_{n=0}^{\infty}\frac{E_{2n+1}^*}{(2n+1)!}2^{2n+1}(-1)^{n+1}.$$
Thus, we have
$$\tan x
=\sum_{n=0}^{\infty}\frac{E_{2n+1}^*}{(2n+1)!}2^{2n+1}(-1)^{n+1}x^{2n+1}.\tag18$$
By (17) and (18), we obtain the following:

\proclaim{ Theorem 3} For $\in \Bbb N$, we have

$$\zeta(2n)=\frac{(-1)^{n-1}(2\pi)^{2n}E_{2n-1}^{*}}{4(2n-1)!(1-4^n)},
$$ where $E_n^*$ are the first kind Euler numbers.
\endproclaim

It is easy to see that
$$\sum_{k=1}^{\infty}\frac{1}{(2k+1)^{2n}}=(1-\frac{1}{4^n})\zeta(2n)
=\frac{(-1)^n(2\pi)^{2n}}{4^{n+1}(2n-1)!}E_{2n-1}^*.$$ Therefore we
obtain the following:

 \proclaim{ Corollary 4} For $n\in\Bbb N $, we
have
$$\sum_{k=1}^{\infty}\frac{1}{(2k+1)^{2n}}=\frac{(-1)^n(2\pi)^{2n}}{4^{n+1}(2n-1)!}E_{2n-1}^*.$$
\endproclaim
Now we try to give the new value of Euler zeta function at positive
integers.
 From the definition of Euler zeta function, we note that
 $$\zeta_E(s)=2\sum_{n=1}^{\infty}\frac{(-1)^n}{n^s}=-2\sum_{n=0}^{\infty}\frac{1}{(2n+1)^s}+\frac{1}{2^{s-1}}
 \zeta(s)
 , \text{ $s\in\Bbb C$ }.\tag19$$
By (19), Theorem 3 and Corollary 4, we obtain the following:

\proclaim{ Theorem 5} For $n\in \Bbb N$, we have
$$\zeta_E(2n)=\frac{(-1)^{n-1}\pi^{2n}(2-4^n)}{2(2n-1)!(1-4^n)}E_{2n-1}^*.$$
\endproclaim
Remark. We note that $\zeta(2)=\frac{\pi^2}{6}$,
$\zeta_{E}(2)=-\frac{\pi^2}{6}$, $\zeta (4)=\frac{\pi^4}{90}$  and
$\zeta_E(4)=-\frac{7\pi^4}{360} \cdots$. For $q\in\Bbb C$ with
$|q|<1,$ $s\in\Bbb C$, $q$-$\zeta$-function is defined by
$$\zeta_q(s)=\sum_{n=1}^{\infty}\frac{q^n}{[n]_q^s}-\frac{1}{s-1}\frac{(1-q)^s}{\log
q}, \text{ cf.[ 12, 17] }  .$$ Note that, $\zeta_q(s)$ is analytic
continuation in $\Bbb C$ with only one simple pole at $s=1$, and
$$\zeta_q(1-k)=-\frac{\beta_{k,q}}{k}, \text{ where $k$ is a
positive integer, cf.[17].}$$ By simple calculation, we easily see
that
$$\aligned
&\sum_{n=1}^{\infty}\frac{(-1)^nq^n}{[n]_q^{2k+1}}\sum_{j=0}^{\infty}\frac{\theta^{2j+1}[n]_q^{2j+1}}{(2j+1)!}
+\frac{1}{\log
q}\sum_{j=0}^{k-1}\frac{(1-q)^{2k-2j}}{(2k-2j-1)(2j+1)!}\\
&-\frac{1}{\log
q}\sum_{j=0}^{k-1}\frac{(1-q)^{2k-2j}}{(2k-2j-1)(2j+1)!}\\
&=\sum_{j=0}^{k-1}\frac{(-1)^j\theta^{2j+1}}{(2j+1)!}\left(-\zeta_q(2k-2j)+\zeta_{q^2}(2k-2j)\frac{2}{[2]_q^{2k-2j}}\right)\\
&-\frac{q}{1+q}\frac{\theta^{2k+1}}{(2k+1)!}(-1)^k+\sum_{j=k+1}^{\infty}\frac{\theta^{2j+1}(-1)^j}{(2j+1)!}
\frac{H_{2j-2k,q}(-q^{-1})}{1+q},
 \endaligned$$
where $H_{n,q}(-q)$ are Carlitz's $q$-Euler numbers with
$\lim_{q\rightarrow 1}H_{n,q}(-q)=E_n^{*},$  cf.[2, 3, 6].
 If $q\rightarrow 1$, then we have
$$\aligned
\sum_{n=1}^{\infty}\frac{(-1)^n}{n^{2k+1}}\sin (n \theta)
&=\sum_{j=0}^{k-1}\frac{(-1)^j}{(2j+1)!}\theta^{2j+1}\left(\frac{2}{2^{2k-2j}}-1\right)\\
&(-1)^{k-j+1}\frac{(2\pi)^{2k-2j}}{2\cdot(2k-2j)!}B_{2k-2j}-\frac{1}{2}\frac{\theta^{2k+1}}{(2k+1)!}(-1)^k.
\endaligned$$
For $k\in\Bbb N$, and $\theta=\frac{\pi}{2}$, it is easy to see
that
$$\sum_{n=1}^{\infty}\frac{(-1)^n}{(2n-1)^{2k+1}}=
\sum_{j=0}^{k-1}
\frac{(-1)^k\pi^{2k+1}(2^{2k-2j}-2)B_{2k-2j}}{(2j+1)!(2k-2j)!2^{2j+2}}-\frac{\pi^{2k+1}(-1)^k}{(2k+1)!2^{2k+2}}.\tag20$$
From (20) and Theorem 1, we can also derive the following equation
(21).
$$\sum_{j=0}^{k-1}
\frac{(-1)^{k-1}\pi^{2k+1}(2^{2k-2j}-2)B_{2k-2j}}{(2j+1)!(2k-2j)!2^{2j+2}}+\frac{\pi^{2k+1}(-1)^k}{(2k+1)!2^{2k+2}}
=(-1)^{k}\frac{E_{2k}}{2(2k)!}\left(\frac{\pi}{2}\right)^{2k+1}.\tag21$$
Thus, we have
$$\sum_{j=0}^{k-1}\frac{(2^{2k-2j}-2)B_{2k-2j}}{(2j+1)!(2k-2j)!2^{2j+2}}=\frac{1}{(2k+1)!2^{2k+2}}-
\frac{E_{2k}}{2^{2k+2}(2k)!}.
$$

\vskip 20pt

\quad Taekyun Kim

\quad EECS, Kyungpook National University,

Taegu 702-701, S. Korea \quad

e-mail:\text{ tkim$\@$knu.ac.kr; tkim64$\@$hanmail.net} \vskip10pt

\enddocument